\magnification=\magstephalf
\parindent=0mm
\parskip=0mm
\hsize=5.75 true in
\baselineskip=15pt
\bigskipamount=24pt plus 6 pt minus 6 pt
\medskipamount=18pt plus 3 pt minus 3 pt
\smallskipamount=12pt plus 2 pt minus 2 pt
\thinmuskip=5mu
\medmuskip=6mu plus 3mu minus 2mu
\thickmuskip=6mu plus 3mu minus 2mu
\overfullrule=0pt
\font\xivrm=cmr17
\def\C{{\mathchoice{\kmplx{6.7pt}}{\kmplx{5.7pt}}{\kmplx{4pt}}{\kmplx{3pt}}}}
\def\kmplx#1{{C\mkern-8.4mu\vrule width 0.4pt height #1\relax\mkern 8mu}}
\def\R{{I\mkern-5mu R}}
\def\N{{I\mkern-5mu N}}
\def\Z{{Z\mkern-8mu Z}}

\def\Any{\,\cdot\,}

\def\norm#1{\mathopen\Vert #1 \mathclose\Vert}

\def\modul#1{\mathopen\vert #1 \mathclose\vert}

\def\loc{\sb{\rm loc}}
\def\Dd#1#2{{d^{#2} \over {d #1}^{#2}}}
\def\square{\mathord{\vbox{\hrule\hbox{\vrule\hskip 9pt\vrule height 9pt}\hrule}}}
\def\litem{\par\noindent\hangindent=\parindent\ltextindent}
\def\ltextindent#1{\hbox to \hangindent{#1\hss}\ignorespaces}
\frenchspacing
\centerline{\xivrm Eigenvalue asymptotics of perturbed periodic}
\centerline{\xivrm Dirac systems in the slow-decay limit.}

\bigskip
\centerline{\it Karl Michael Schmidt}
\centerline{\it School of Mathematics, Cardiff University, 23 Senghennydd Rd}
\centerline{\it Cardiff CF24 4YH, UK}

\footnote{}{{\bf 2000 MSC:} 34L20, 34L40, 47E05, 81Q10, 81Q15}

\bigskip
{\bf Abstract.}
A perturbation decaying to $0$ at $\infty$ and not too irregular at $0$
introduces at most a discrete set of eigenvalues into the spectral gaps of
a one-dimensional Dirac operator on the half-line.
We show that the number of these eigenvalues in a compact subset of a
gap in the essential spectrum is given by a quasi-semiclassical asymptotic
formula in the slow-decay limit, which for power-decaying perturbations
is equivalent to the large-coupling limit.
This asymptotic behaviour elucidates the origin of the dense point spectrum
observed in spherically symmetric, radially periodic three-dimensional
Dirac operators.

\bigskip

{\bf 1 Introduction.}

\bigskip
For a large class of potentials, the semiclassical Weyl formula gives a
correct asymptotic description of the total multiplicity of the lower
spectrum of a Schr\"odinger operator in the large coupling limit (see
[18] Chapter XIII.15 and the references given there).
However, there are also some notable exceptions, e.g. in the
two-dimensional case ([9], [10]).
In recent years an analogous asymptotic analysis of the point spectrum
arising in spectral gaps of Schr\"odinger operators under perturbations
has attracted considerable attention; starting from [1], Birman has
developed a general framework to study this problem [2], [3], [4], [5],
[6], [7]; see also [16].

More specifically, Sobolev [26] has studied the perturbed periodic
one-dimensional Schr\"o\-din\-ger operator, showing that for a wide range
of power-decaying perturbations, the number of eigenvalues arising in a
closed subinterval of a spectral gap of the unperturbed problem is
asymptotically given by a quasi-semiclassical formula in which the
quasimomentum of the periodic background problem takes the role of the
ordinary momentum in the Weyl formula.
This result has an interesting application to spherically symmetric
Schr\"odinger operators in $\R^n$ with a radially periodic potential,
which have dense point spectrum in all spectral gaps of the corresponding
one-dimensional periodic operator [14].
Indeed, treating the angular momentum term of the one-dimensional
half-line operators in the partial-wave decomposition
$$
  - \Delta + q(\modul\Any)
  = \bigoplus_{l \in \N_0}
  - \Dd{r}{2} + q(r) + \big(l (l + n - 2) + {(n - 1) (n - 3) \over 4}\big)
      /r^2
\qquad (r \in (0, \infty)),
$$
with periodic $q$, as a perturbation, the quasi-semiclassical formula explains
the origin of the dense eigenvalues, and describes their asymptotic
density in the limit $l \rightarrow \infty$.
Furthermore, numerical experiments have shown that the asymptotic
formula can give a surprisingly accurate prediction for the number of
eigenvalues even for small values of the coupling constant [11].
(We mention in passing that the angular momentum term also happens to be
a boundary case for the question whether the total number
of eigenvalues introduced in a spectral gap for an individual half-line
operator is finite or infinite ([19], [20], [21], [22], [25]); but this
critical behaviour is not apparent in the large-coupling limit.)

The relativistic counterpart of the Schr\"odinger operator, the Dirac
operator, is unbounded below, so that there is no lower spectrum,
and one is always in a gap situation when studying the discrete spectrum,
even when there is no periodic background potential (for the large coupling
asymptotics in that case, see e.g. [8], who correct a result by [15],
and [12]).
Again, an interesting class of perturbed periodic one-dimensional operators
arises from the partial-wave decomposition of the spherically symmetric,
radially periodic Dirac operator in $\R^3$,
$$
  H =
  -i \alpha \cdot \nabla + m \beta + q(\modul\Any)
  \cong \bigoplus_{k \in \Z \setminus\{0\}}
  -i \sigma_2 \Dd{r}{} + m \sigma_3 + q(r) + {k \over r} \sigma_1
$$
where $\alpha_1, \alpha_2, \alpha_3, \beta$ are $4\times 4$ Dirac matrices,
and
$$
  \sigma_1 = \pmatrix{0 & 1 \cr 1 & 0 \cr},
\quad
  \sigma_2 = \pmatrix{0 & -i \cr i & 0 \cr},
\quad
  \sigma_3 = \pmatrix{1 & 0 \cr 0 & -1 \cr}.
$$
This operator also has dense point spectrum in the spectral gaps of the
corresponding one-dimensional periodic operator, but one exceptional
gap may or may not contain dense point spectrum --- it is an open
question whether partial filling of the gap can occur [23], [24].
The underlying reason for the complications in the Dirac case is the fact
that the angular momentum term $\sigma_1 k/r$, in contrast to the
Schr\"odinger angular momentum term, does not have a definite sign.
Since it does not depend monotonically on the quantum number $k$, the
eigenvalues can move either way as $k$ increases.
Methods such as Sturm comparison, which yields a fairly quick
proof for the quasi-semiclassical formula for perturbed Sturm-Liouville
operators, are of little use when applied to matrix perturbations
of Dirac operators.

In the present paper we study the asymptotic distribution of eigenvalues
in spectral gaps of periodic Dirac operators with a matrix perturbation of
the type of the angular momentum term.
We assume that the perturbation tends to $0$ at $\infty$, but no assumptions
on its sign or rate of decay are needed.
The asymptotic limit we consider is in fact a `slow-decay' limit, which for
perturbations of inverse power decay is equivalent to the `large-coupling'
limit, an observation central to Sobolev's result.
More explicitly, given a Dirac system
$$
  \tau = - i \sigma_2 \Dd{x}{} + m \sigma_3 + q(x)
$$
with $m > 0$ and real-valued, $\alpha$-periodic $q \in L^1\loc(\R)$, and
a real-valued perturbation template
$$
  l_0 \in C(0, \infty),
\quad
  \lim_{\varrho \rightarrow 0} l_0(\varrho) = \infty,
\quad
  \lim_{\varrho \rightarrow \infty} l_0(\varrho) = 0,
$$
we investigate the number of eigenvalues of the perturbed Dirac
system on $(0, \infty)$
$$
  \tau(c) = -i \sigma_2 \Dd{r}{} + m \sigma_3 + q(r) + l_0(r/c) \sigma_1
$$
in a compact subinterval of a spectral gap asymptotically
as $c \rightarrow \infty$.

The angular momentum term is clearly of this type, with
$l_0(\varrho) = 1/\varrho$.

We follow Sobolev's basic idea of comparing $\tau(c)$ to an operator
in which $l_0(r/c)$ is locally replaced by a constant, but in the
absence of Sturm comparison to obtain upper and lower bounds we must
resort to the coarser instrument of operator perturbation theory.
This complication leads to severe difficulties near $0$, where the
perturbation is divergent: in the Sturm-Liouville case this actually
helps, as the perturbed eigenvalue equation becomes disconjugate near $0$,
so that no spectrum is produced near that end-point.
The case of the Dirac system, however, with an unperturbed
operator unbounded below and a strong perturbation of no definite sign,
is entirely different.
For this reason we introduce the following additional assumptions on $q$
and $l_0$ to ensure that essentially no spectrum is created near $0$; as
will be seen, they also imply that $\tau(c)$ is in the limit point case at
$0$ for sufficiently large $c$.
$$
  q \in L^\infty(\R);
\qquad
  l_0 \in AC\loc(0, \Any),
\quad
  \limsup_{\varrho \rightarrow 0}
  {\modul{l_0'(\varrho)} \over l_0^2(\varrho)} < \infty.
\leqno (H)
$$
The requirements on $l_0$ are clearly satisfied in the case of the angular
momentum term.

Let $k(\lambda, l)$ be the quasi-momentum (cf. Section 2 below) of the
periodic equation
$$
  \tau u = (\lambda - l \sigma_1) u,
$$
for $\lambda, l \in \R$.

Our main result is

\medskip
{\bf Theorem 1.}\qquad
{\it
Assume that {\/\rm (H)} holds in addition to the general hypotheses, and
let $[\lambda_1, \lambda_2]$ be compactly contained in a spectral gap of
the self-adjoint realisation $h$ of $\tau$.

Then for sufficiently large $c > 0$, $\tau(c)$ is essentially self-adjoint
on $C_0^\infty(0, \infty)$, and the number $N(c)$
of eigenvalues in $[\lambda_1, \lambda_2]$ of its self-adjoint extension
$h(c)$ satisfies
$$
  \lim_{c \rightarrow \infty} {N(c) \over c}
  = {1 \over \alpha \pi} \int_0^\infty (k(\lambda_2, l_0(\varrho)) -
         k(\lambda_1, l_0(\varrho)))\,d\varrho.
$$
\/}

\smallskip
{\it Remarks.\/}

1.
Our assumptions imply that the essential spectra of $h(c)$ and
$h$ coincide, so $h(c)$ only has discrete eigenvalues in $[\lambda_1,
\lambda_2]$.
It will be apparent from the proof that the asymptotic formula in Theorem 1
continues to hold if, instead of (H) and $\lim\limits_{\varrho\rightarrow 0}
l_0(\varrho) = \infty$, we assume $l_0 \in C([0, \infty))$.

2.
In [23] Theorem 2, it was shown that $\lambda \in \R$ is a point of the
essential spectrum of $H$ if $\lambda$ is in the spectrum of the periodic
operator $h + l \sigma_1$ for some $\l \in \R$, i.e. if $\lambda$ is a
point of growth of $k(\Any, l)$.
In a rather weak sense, Theorem 1 is a reverse of this statement:
if a $\lambda$ interval does not intersect the spectrum of $h + l \sigma_1$
for any $l \in \R$, then the integral on the r.h.s. of the asymptotic formula
vanishes,
i.e. the asymptotic eigenvalue density is $0$.
Of course, this does not rule out the existence of dense point spectrum
of $H$ in this interval.

\bigskip
{\bf 2 Rotation number and quasimomentum of periodic Dirac systems.}

\bigskip
The qualitative behaviour of solutions of a linear ordinary differential
equation with periodic coefficients is characterised by its monodromy matrix,
the value of the canonical fundamental system after one period.
In the case of an $\alpha$-periodic Sturm-Liouville or Dirac equation,
this is a $2 \times 2$ matrix with real entries and (Wronskian) determinant 1;
the position of its eigenvalues $\mu_1, \mu_2$ in the complex plane
is thus fully determined by its trace, the {\it discriminant\/} $D =
\mu_1 + \mu_2$. If $\mu_1 \neq \mu_2$, then there are corresponding Floquet
solutions $u_1, u_2$ satisfying
$$
  u_j(x + \alpha) = \mu_j\,u_j(x)
\qquad
  (x \in \R, j \in \{1, 2\}).
$$
Hence it is easy to see that the equation is stable (all solutions are
globally bounded) if $\modul{D} < 2$, and unstable (all solutions are
unbounded) if $\modul{D} > 2$; as $D$ is an analytic function of the
spectral parameter, the real line splits into alternating stability and
instability intervals.
This already provides a complete description of the (purely absolutely
continuous) spectrum of the corresponding self-adjoint ordinary differential
operator: it is the closure of the union of all stability intervals,
spectral gaps corresponding to non-degenerate instability intervals
(cf. [13], [27]).

Furthermore, the oscillation behaviour of solutions is very closely linked
to the discriminant as well.
Thus for the periodic Dirac system, it is not difficult to verify along the
lines of [13] proof of Thm 3.1.2, that for any value of the spectral
parameter in the closure of the $n$th instability interval, $n \in \Z$,
the Pr\"ufer angle $\vartheta$ of any $\R^2$-valued solution $u$,
defined by
$$
  u = R \pmatrix{\cos\vartheta \cr -\sin\vartheta \cr}
$$
with $R > 0$, satisfies
$$
  \vartheta(x) = {n \pi x \over \alpha} + O(1)
\qquad
  (x \rightarrow \infty).
$$
In particular, the asymptotic rate of growth of $\vartheta$, the so-called {\it
rotation number\/}
$$
 \lim_{x\rightarrow \infty} {\vartheta(x) \over x} = {n \pi \over \alpha},
$$
is constant in instability intervals.

It turns out that the rotation number is in fact well-defined as a
continuous non-decreasing function of the spectral parameter on the whole
real line, which can be expressed in terms of the discriminant in the
stability intervals.
For periodic Sturm-Liouville equations this was elegantly shown in [17],
using a connection between the winding number of a complex-valued
Floquet solution in the punctured plane and the Pr\"ufer
angle of its real part.
Unfortunately, this argument does not carry over to the Dirac system,
with Floquet solutions (in the case of stability) consisting of two
complex-valued components with no immediate link between them.
Nevertheless, the following statement holds.

\medskip
{\bf Theorem 2.}\qquad
{\it
Let $\vartheta$ be the Pr\"ufer angle of an $\R^2$-valued solution of the
Dirac system
$$
  (-i \sigma_2 \Dd{x}{} + m \sigma_3 + l \sigma_1 + q) u = \lambda u
$$
with real-valued, $\alpha$-periodic coefficients $m, l, q \in L^1\loc(\R)$,
$\lambda \in \R$ in a stability interval, and let $D$ be the discriminant.

Then there is $k(\lambda) \in \R$ such that
$\vartheta(x) = {k(\lambda) x \over \alpha} + O(1)$
$(x \rightarrow \infty)$, and $D = 2 \cos k(\lambda)$.
\/}

\medskip
$k(\lambda)$ is called the {\it quasimomentum\/} of this periodic equation;
$k(\lambda) / \alpha$ is the rotation number.
$k(\lambda) / (\alpha \pi)$ is called the {\it integrated density of states\/}
in view of the following consequence of Theorem 2, which is readily
obtained using oscillation theory ([27] Thm 14.7, 14.8).

\medskip
{\bf Corollary 1.}\qquad
{\it
In the situation of Theorem 2,
let $(a_n)_{n\in\N}$, $(b_n)_{n \in \N}$ be sequences of real numbers such
that $a_n < b_n$ and $\lim\limits_{n \rightarrow \infty} b_n - a_n = \infty$,
and $t_n$ any self-adjoint realisation of
$$
  -i \sigma_2 \Dd{x}{} + m \sigma_3 + l \sigma_1 + q
$$
on $[a_n, b_n]$ with separated boundary conditions.
Then for $\lambda_1 < \lambda_2$ the
number $N_n$ of eigenvalues of $t_n$ in $[\lambda_1, \lambda_2]$ satisfies
$$
  \lim_{n\rightarrow\infty} {N_n \over b_n - a_n}
  = {k(\lambda_2) - k(\lambda_1) \over \alpha \pi}.
$$
\/}

\medskip
In the proof of Theorem 2 we use the following elementary observation on
complex solutions of general Dirac systems with real-valued coefficients.

\medskip
{\bf Lemma 1.}\qquad
{\it
Let $u : I \rightarrow \C^2$ be a non-trivial solution of
$$
  (-i \sigma_2 \Dd{x}{} + m \sigma_3 + l \sigma_1 + q) u = \lambda u
$$
($m, l, q \in L^1\loc(I)$ real-valued, $\lambda \in \R$).
Then the following statements are equivalent:
\parindent=.2in

\litem{a)}
There is $x \in I$, $\varphi \in C$, $\modul\varphi = 1$, such that
$\varphi u(x) \in \R^2$.

\litem{b)}
There is $\varphi \in \C$, $\modul\varphi = 1$, such that
$\varphi u : I \rightarrow \R^2$.

\litem{c)}
$u, \overline{u}$ are linearly dependent.
\/}

\bigskip
{\it Proof\/} of Theorem 2.\qquad

As $\modul{D} < 2$, the monodromy matrix has complex eigenvalues
$\mu \neq \overline{\mu}$, $\modul{\mu} = 1$.
By Lemma 1, the corresponding Floquet solutions $u$, $\overline{u}$
are linearly independent as $\mu \notin \R$.

Again by Lemma 1, the components $u_1, u_2$ have no zeros, and
their arguments (locally absolutely continuous functions) are
nowhere equal mod $\pi$. By choosing an appropriate branch of the
argument, we can assume $0 < \modul{\arg u_1 - \arg u_2} < \pi$
throughout. More specifically, there are locally absolutely
continuous functions $R_1, R_2 > 0$, $\varphi_1, \varphi_2$, and
$\nu \in \{-1, 1\}$ such that $u_1 = R_1 e^{i\varphi_1}$, $u_2 =
i \nu R_2 e^{i \varphi_2}$, and $\modul{\varphi_1 - \varphi_2} <
\pi/2$.
Interchanging $u$ and $\overline u$ if necessary, we can assume
without loss of generality that $\nu = 1$.

Consider the $\R^2$-valued solution $v := \mathop{\rm Re} u$; its Pr\"ufer
angle satisfies
$$
  \tan \vartheta = -{v_2 \over v_1}
  = {R_2 \sin \varphi_2 \over R_1 \cos \varphi_1}
  = {R_2 \over R_1}\,\cos (\varphi_2 - \varphi_1)\,
    (\tan \varphi_1 + \tan (\varphi_2 - \varphi_1)).
$$
Thus $\vartheta$ and $\varphi_1$ are connected by a modified Kepler
transformation (cf. [25]); as a consequence, they take the values
$(\Z + 1/2) \pi$ at the same points, and their difference is globally
bounded.

From $u_1(x + \alpha) = \mu u_1(x)$ we find that $R_1$ is $\alpha$-periodic,
and that there is $k \in \R$ such that $\mu = e^{ik}$ and
$\varphi_1(x + \alpha) = \varphi_1(x) + k$ $(x \in \R)$.
Thus $\varphi_1(x) - {kx \over \alpha}$ is $\alpha$-periodic and continuous,
hence globally bounded, and it follows that
$$
  \vartheta(x) = {xk \over \alpha} + O(1)
\qquad
  (x \rightarrow \infty).
$$
Furthermore, $D = \mu + \overline{\mu} = 2 \cos k$.
\hfill $\square$

\bigskip
{\bf 3 Proof of Theorem 1.}

\nobreak
\bigskip
We shall use the following consequences of the spectral theorem
(cf. [24] Lemma 6).

\medskip
{\bf Lemma 2.}\qquad
{\it
Let $L$ be a self-adjoint operator with purely discrete spectrum of finite
total multiplicity $N_I(L)$ in the real interval
$I = [\lambda_1, \lambda_2]$.

\parindent=1em

\smallskip
\litem{a)} {\rm(Decomposition Principle)}
If $L$ is a $k$-dimensional extension of a restriction of a self-adjoint
operator $L'$, then $N_I(L) \le N_I(L') + k$.

\smallskip
\litem{b)} {\rm(Bounded Perturbations)}
If $A$ is symmetric and bounded, then $N_I(L) \le N_{I'}(L + A)$,
where $I' = [\lambda_1 - \norm{A}, \lambda_2 + \norm{A}]$.

}

\medskip
The proof of Theorem 1 proceeds as follows.
We introduce the rescaled independent variable
$\varrho = r/c \in (0, \infty)$,
and use the decomposition principle to split off the two regions near the
singular end-points $0$ and $\infty$.
The regular operator on the intermediate interval is then, by another
application of the decomposition principle, split into operators on a
large but finite number of subintervals on which the perturbation
$\sigma_1 l_0(\varrho)$ does not change very much.
Lemma 2 b) relates the number of eigenvalues
of each subinterval operator to that of an operator with periodic
coefficients on the same subinterval.
By virtue of Corollary 1, the asymptotic number of eigenvalues
for the latter can be expressed by means of the quasimomentum,
considering that in terms of the original variable $r$, the subintervals
grow beyond all bounds in the limit $c \rightarrow \infty$.

The operator on the unbounded interval near $\infty$, on which the
perturbation is small, is easily dealt with in a similar manner.
The end-point $0$, however, presents more serious difficulties, which
reflect the lack of monotonicity of spectral behaviour with respect to a
matrix perturbation; and it is at this point only that our rather strong
hypotheses (H) enter.
We show that with a suitably chosen boundary
condition, the operator on the interval near $0$ has no eigenvalues at all
in $[\lambda_1, \lambda_2]$ for sufficiently large $c$.
The asymptotic formula then follows in the limit of infinite refinement
of the subintervals.

\medskip
The Dirac system
$$
  \tau(c) u = \lambda u
$$
is in the limit point case at infinity ([27] Cor. to Theorem 6.8).
We begin the {\it proof\/} of Theorem 1 by showing that for sufficiently
large $c$, it is in the limit point case at $0$ as well,
thus proving the essential self-adjointness.

Let $\hat\varrho > 0$ and $C > 0$ be such that $l_0(\varrho) > 0$
$(\varrho \in (0, \hat\varrho))$, and
$$
  {\modul{l_0'(\varrho)} \over l_0(\varrho)^2} \le C
\qquad
  (\varrho \in (0, \hat\varrho)).
$$
Then $(-\log l_0)'(\varrho) \le C l_0(\varrho)$, and hence
$$
  l_0(\varrho)
  \le l_0(\hat\varrho) \exp\left(C\int_\varrho^{\hat\varrho} l_0\right)
\qquad
  (\varrho \in (0, \hat\varrho));
$$
and similarly
$$
  l_0(\varrho) - l_0(\hat\varrho) = - \int_\varrho^{\hat\varrho} l_0'
  \le C\int_\varrho^{\hat\varrho} l_0^2,
$$
which implies $l_0 \notin L^2(0, \Any)$ in view of
$\lim\limits_{\varrho\rightarrow 0} l_0(\varrho) = \infty$.

As $m \sigma_3 + q - \lambda$ is bounded, it is sufficient to establish that
$$
  (-i \sigma_2 \Dd{r}{} + \sigma_1 l_0(r/c)) u = 0
$$
is in the limit point case at $0$.
If $u$ is a non-trivial $\R^2$-valued solution of this equation, then
$u_1'(r) = - l_0(r/c) u_1(r)$, and setting $v(\varrho) := u_1(c \varrho)$,
we find for $c \ge C$
$$
  v(\varrho)
  = v(\hat\varrho) \exp\left(c \int_\varrho^{\hat\varrho} l_0 \right)
  \ge {v(\hat\varrho) \over l_0(\hat\varrho)} l_0(\varrho),
$$
so $v \notin L^2(0, \Any)$.

To prove the asymptotic formula, we first observe that the region near $0$
does not contribute to the eigenvalue count.
Indeed, let $c_0 := C+1$; then
$$
  l_0^2(\varrho) - {1 \over c_0} \modul{l_0'(\varrho)}
  \ge {l_0^2(\varrho) \over c_0} \rightarrow \infty
\qquad
  (\varrho \rightarrow 0),
$$
and thus there is $\varrho_0 > 0$ such that
$$
  l_0^2(\varrho) - {1 \over c_0} \modul{l_0'(\varrho)}
  \ge (\norm{q}_\infty + \max\{\modul{\lambda_1}, \modul{\lambda_2}\} + 1)^2.
\qquad
  (\varrho \in (0, \varrho_0)).
$$

Let $c \ge c_0$, and $h_0(c)$ the self-adjoint realisation of $\tau(c)$ on
$(0, c \varrho_0)$ with the boundary condition
$$
  u_1(c \varrho_0) + u_2(c \varrho_0) = 0.
$$
As $\tau(c)$ is in the limit point case at $0$, $h_0(c)$ is essentially
self-adjoint
on $D_0 := \{u \in D(h_0(c)) \mathrel{|} u \equiv 0 \hbox{ near } 0\}$.
Hence if $\lambda \in \sigma(h_0(c))$ and $\varepsilon \in (0, 1)$, there is
$u \in D_0 \setminus \{0\}$ such that
$\norm{(h_0(c) - \lambda) u} < \varepsilon \norm{u}$,
so
$$
  \norm{-i \sigma_2 u' + (m \sigma_3 + l_0(\Any / c) \sigma_1) u}
  \le (\norm{q}_\infty + \modul{\lambda} + \varepsilon) \norm{u}.
$$
But on the other hand, integration by parts yields
$$\eqalign{
  &\norm{-i \sigma_2 u' + (m \sigma_3 + l_0(\Any / c) \sigma_1) u}^2
\cr
  &\qquad
  =\int_0^{c \varrho_0} (\modul{u'}^2 + (m^2 + l_0(\Any/c)^2) \modul{u}^2
    - m ((\sigma_1 u)^T \overline{u})'
    + l_0(\Any/c) ((\sigma_3 u)^T \overline{u})')
\cr
  &\qquad
  = m \modul{u(c \varrho_0)}^2
  +\int_0^{c \varrho_0} (\modul{u'}^2 + (m^2 + l_0(\Any/c)^2) \modul{u}^2
    - {1 \over c} l_0'(\Any/c) ((\sigma_3 u)^T \overline{u}))
\cr
  &\qquad
  \ge \int_0^{c \varrho_0} (l_0(r/c)^2 - {1 \over c} \modul{l_0'(r/c)})
    \modul{u(r)}^2\,dr
  \ge (\norm{q}_\infty + \max\{\modul{\lambda_1}, \modul{\lambda_2}\} + 1)^2
    \norm{u}^2.
\/}$$
Thus $\modul\lambda > \max \{\modul{\lambda_1}, \modul{\lambda_2}\}$; in particular,
$h(c)$ has no spectrum in $[\lambda_1, \lambda_2]$.

A similar calculation shows that for all fixed $\varrho \in (0, \varrho_0)$,
the self-adjoint periodic operator
$$
  -i \sigma_2 \Dd{x}{} + m \sigma_3 + q(x) + l_0(\varrho) \sigma_1
$$
has a spectral gap, and hence instability interval, containing $[\lambda_1,
\lambda_2]$, and it follows that
$k(\lambda_1, l_0(\varrho)) = k(\lambda_2, l_0(\varrho))$ $(\varrho \in (0, \varrho_0))$.

\medskip
We now turn to the remaining interval $(c \varrho_0, \infty)$.
Choose $\delta_\infty \in (0, 1)$ so small that
$\delta_\infty < (\lambda_2 - \lambda_1)/2$, and
$[\lambda_1 - \delta_\infty, \lambda_2 + \delta_\infty]$ is still compactly
contained in the spectral gap of the periodic problem.
Fix $P_0 > \varrho_0$ such that $\modul{l_0(\varrho)} < \delta_\infty$
$(\varrho > P_0)$.

Let $\varepsilon > 0$. There is $L > 0$ such that $\modul{l_0(\varrho)}
\le L$ $(\varrho \in (\varrho_0, \infty))$.
As the quasimomentum $k(\lambda, l)$ is uniformly continuous
on the compact set $K := [\lambda_1 - 1, \lambda_2 + 1] \times [-L, L]$,
there is $\delta \in (0, \delta_\infty)$ such that
$$
  (\mu, l), (\mu', l') \in K, \modul{\mu - \mu'}, \modul{l - l'} < \delta
  \mathrel{\Rightarrow} \modul{k(\mu, l) - k(\mu', l')} < \varepsilon.
$$
As $l_0$ is uniformly continuous on $[\varrho_0, P_0]$, there is $\gamma > 0$
such that
$$
  \modul{l_0(x) - l_0(y)} < \delta
\qquad
  (x, y \in [\varrho_0, P_0], \modul{x - y} < \gamma).
$$

Now consider a partitioning of the interval $I := (\varrho_0, P_0)$ into
$n$ subintervals $I_j = (\varrho_{j-1}, \varrho_j)$ with
$\modul{I_j} < \gamma$ $(j \in \{1, \dots, n\})$.

Let $h_j(c)$, $j \in \{1, \dots, n\}$, and $h_\infty(c)$ be self-adjoint
realisations of $\tau(c)$ on $c I_j$, and $(c P_0, \infty)$, respectively,
with separated boundary conditions at the regular end-points.
Then $h_0(c) \oplus \bigoplus\limits_{j=1}^n h_j(c) \oplus h_\infty(c)$
is a $2(n+1)$-dimensional extension of the $2(n+1)$-dimensional restriction
of $h(c)$ with domain
$\{u \in D(h(c)) \mathrel{|} u(c \varrho_j) = 0 (j \in \{0, \dots, n\})\}$;
comparing these operators by means of the decomposition
principle (Lemma 2 a) we find
$$
  \modul{N_0(c) + \sum_{j=1}^n N_j(c) + N_\infty(c) - N(c)} \le 2(n+1)
$$
for the total spectral multiplicities in $[\lambda_1, \lambda_2]$.
As observed above, $N_0(c) = 0$.

To estimate $N_j(c)$ for $j \in \{1, \dots, n\}$, choose $\tilde\varrho_j \in I_j$,
let $l_j := l_0(\tilde\varrho_j)$, and define self-adjoint operators
$$
  \tilde h_j(c) = -i \sigma \Dd{r}{} + m \sigma_3 + q(r) + l_j \sigma_1
$$
on $c I_j$ with $D(\tilde h_j(c)) = D(h_j(c))$.
Then Lemma 2 b) and $\modul{l_0(\varrho) - l_j} < \delta$
$(\varrho \in I_j)$ yield the bounds
$$
  \tilde N_j(\lambda_1 + \delta, \lambda_2 - \delta; c)
  \le N_j(c) \le
  \tilde N_j(\lambda_1 - \delta, \lambda_2 + \delta; c),
$$
where $\tilde N_j(\mu_1, \mu_2; c)$ is the number of eigenvalues of
$\tilde h_j(c)$ in $[\mu_1, \mu_2]$.

As $\tilde h_j(c)$ has $\alpha$-periodic coefficients, Corollary 1 implies
$$
  \lim_{c \rightarrow \infty}
    {\tilde N_j(\mu_1, \mu_2; c) \over c \modul{I_j}}
  = {k(\mu_2, l_j) - k(\mu_1, l_j) \over \alpha \pi}.
$$
Thus we find
$$\eqalign{
  &{\modul{I_j} \over \alpha \pi}
    (k(\lambda_2, l_j) - k(\lambda_1, l_j) - 2 \varepsilon)
  \le
  {\modul{I_j} \over \alpha \pi}
    (k(\lambda_2 - \delta, l_j) - k(\lambda_1 + \delta, l_j))
\cr
  &\qquad
  = \lim_{c \rightarrow \infty} {1 \over c}
     \tilde N_j(\lambda_1 + \delta, \lambda_2 - \delta; c)
  \le \liminf_{c \rightarrow \infty} {N_j(c) \over c}
\cr
  &\qquad
  \le \limsup_{c \rightarrow \infty} {N_j(c) \over c}
  \le \lim_{c \rightarrow \infty} {1 \over c}
     \tilde N_j(\lambda_1 - \delta, \lambda_2 + \delta; c)
\cr
  &\qquad
  =
  {\modul{I_j} \over \alpha \pi}
    (k(\lambda_2 + \delta, l_j) - k(\lambda_1 - \delta, l_j))
  \le
  {\modul{I_j} \over \alpha \pi}
    (k(\lambda_2, l_j) - k(\lambda_1, l_j) + 2 \varepsilon).
\cr}$$

In order to estimate $N_\infty(c)$, we observe that
$$
  \tilde N_\infty(\lambda_1 + \delta_\infty, \lambda_2 - \delta_\infty; c)
  \le N_\infty(c) \le
  \tilde N_\infty(\lambda_1 - \delta_\infty, \lambda_2 + \delta_\infty; c),
$$
where $\tilde N_j(\mu_1, \mu_2; c)$ is the number of eigenvalues in
$[\mu_1, \mu_2]$ of the self-adjoint realisation
$\tilde h_\infty(c)$ of $\tau$ on $(c P_0, \infty)$ with
$D(\tilde h_\infty(c)) = D(h_\infty(c))$.
Since $\tilde h_\infty(c)$ is a direct summand of a 2-dimensional
extension of a restriction of the unperturbed periodic operator on $\R$,
we find, again by the decomposition principle, that
$0 \le \tilde N_\infty(\mu_1, \mu_2; c) \le 2$ for any $[\mu_1, \mu_2]$
compactly contained in a spectral gap.

Summing up, we obtain
$$\eqalign{
  &{1 \over \alpha\pi} \sum_{j=1}^n (k(\lambda_2, l_j) - k(\lambda_1, l_j))
   \modul{I_j} - {2 \varepsilon \over \alpha\pi}\modul{I}
  \le \liminf_{c \rightarrow \infty} {1 \over c} \sum_{j=1}^n N_j(c)
\cr
  &\qquad
  = \liminf_{c \rightarrow \infty} {N(c) \over c}
  \le \limsup_{c \rightarrow \infty} {N(c) \over c}
  = \limsup_{c \rightarrow \infty} {1 \over c} \sum_{j=1}^n N_j(c)
\cr
  &\qquad
  \le
  {1 \over \alpha\pi} \sum_{j=1}^n (k(\lambda_2, l_j) - k(\lambda_1, l_j))
   \modul{I_j} + {2 \varepsilon \over \alpha\pi}\modul{I},
\cr}$$

and refining the Riemann sums,
$$\eqalign{
  &{1 \over \alpha\pi}
   \int_I (k(\lambda_2, l_0(\varrho)) - k(\lambda_1, l_0(\varrho)))\,d\varrho
  - {2 \varepsilon \over \alpha \pi} \modul{I}
  \le \liminf_{c \rightarrow \infty} {N(c) \over c}
\cr
  &\qquad
  \le \limsup_{c \rightarrow \infty} {N(c) \over c}
  \le {1 \over \alpha\pi}
   \int_I (k(\lambda_2, l_0(\varrho)) - k(\lambda_1, l_0(\varrho)))\,d\varrho
  + {2 \varepsilon \over \alpha \pi} \modul{I}.
\cr}$$
The assertion of Theorem  1 follows, as $\varepsilon > 0$ is arbitrary, and
$k(\lambda_2, l_0(\varrho)) = k(\lambda_1, l_0(\varrho))$
$(\varrho \in (0, \infty) \setminus I)$.

This concludes the proof of Theorem 1.

\bigskip
{\bf References.}

\parindent=.3 in

\litem{1}
Alama S., Deift P.A., Hempel R. Eigenvalue branches of the Schr\"odinger
operator $H - \lambda W$ in a gap of $\sigma(H)$. {\it Commun. Math. Phys.\/}
{\bf 121} (1989) 291--321

\litem{2}
Birman M.Sh. Discrete spectrum in the gaps of the continuous one in the
large-coupling-constant limit. in: {\it Order, disorder and chaos in quantum
systems (Dubna 1989), Oper. Theory Adv. Appl.\/} {\bf 46}, Birkh\"auser,
Basel 1990, pp. 17--25

\litem{3}
Birman M.Sh. Discrete spectrum in the gaps of a continuous one for
perturbations with large coupling limit. {\it Adv. Sov. Math.\/} {\bf 7}
(1991) 57--73

\litem{4}
Birman M.Sh. On a discrete spectrum in gaps of a second order perturbed
periodic operator. {\it Funct. Anal. Appl.\/} {\bf 25} (2) (1991) 158--161

\litem{5}
Birman M.Sh. The discrete spectrum in gaps of the perturbed periodic
Schr\"odinger operator I. Regular perturbations. in: {\it Boundary value
problems, Schr\"odinger operators, deformation quantization. Math. Top.\/}
{\bf 8}, Akademie Verlag, Berlin 1995, pp. 334--352

\litem{6}
Birman M.Sh. The discrete spectrum of the periodic Schr\"odinger operator
perturbed by a decreasing potential. {\it St. Petersburg Math. J.\/}
{\bf 8} (1) (1997) 1--14

\litem{7}
Birman M.Sh. The discrete spectrum in gaps of the perturbed periodic
Schr\"odinger operator II. Nonregular perturbations. {\it St. Petersburg
Math. J.\/} {\bf 9} (6) (1998) 1073--1095

\litem{8}
Birman M.Sh., Laptev A. Discrete spectrum of the perturbed Dirac operator.
{\it Ark. Mat.\/} {\bf 32} (1994) 13--32

\litem{9}
Birman M.Sh., Laptev A. The negative discrete spectrum of a two-dimensional
Schr\"odinger operator. {\it Comm. Pure Appl. Math.\/} {\bf 49} (1996)
967--997

\litem{10}
Birman M.Sh., Laptev A., Solomyak M. The negative discrete spectrum of the
operator $-\Delta^l - \alpha V$ in $L_2(\R^d)$ for $d$ even and $2l \ge d$.
{\it Ark. Mat.\/} {\bf 35} (1997) 87--126

\litem{11}
Brown B.M., Eastham M.S.P., Hinz A.M., Schmidt K.M. Distribution of
eigenvalues in gaps of the essential spectrum of Sturm-Liouville operators
--- a numerical approach. {\it J. Comp. Anal. Appl.\/} (to appear)

\litem{12}
Cancelier C., L\'evy-Bruhl P., Nourrigat J. Remarks on the spectrum of
Dirac operators. {\it Acta Appl. Math.\/} {\bf 45} (1996) 349--364

\litem{13}
Eastham M.S.P. {\it The spectral theory of periodic differential equations.\/}
Scottish Academic Press, Edinburgh 1973

\litem{14}
Hempel R. Herbst I., Hinz A.M., Kalf H. Intervals of dense point spectrum for
spherically symmetric Schr\"odinger operators of the type
$-\Delta + \cos \modul{x}$. {\it J. London Math. Soc. (2)\/} {\bf 43} (1989)
295--304

\litem{15}
Klaus M. On the point spectrum of Dirac operators. {\it Helv. Phys. Acta\/}
{\bf 53} (1980) 453--462

\litem{16}
Laptev A. Asymptotics of the negative discrete spectrum of a class of
Schr\"odinger operators with large coupling constant. {\it Proc. Amer.
Math. Soc. (2)\/} {\bf 119} (1993) 481--488

\litem{17}
Moser J. An example of a Schr\"odinger equation with almost periodic
potential and nowhere dense spectrum. {\it Comment. Math. Helv.\/} {\bf 56}
(1981) 198--224

\litem{18}
Reed M, Simon B. {\it Methods of modern mathematical physics IV:
Analysis of operators.\/} Academic Press, New York 1978

\litem{19}
Rofe-Beketov F.S. Spectral analysis of the Hill operator and of its
perturbations. {\it Functional analysis\/} {\bf 9} (1977) 144-155 (Russian)

\litem{20}
Rofe-Beketov F.S. A generalisation of the Pr\"ufer transformation and
the discrete spectrum in gaps of the continuous spectrum. in: {\it Spectral
theory of operators.\/} Elm, Baku 1979, pp. 146--153 (Russian)

\litem{21}
Rofe-Beketov F.S. Spectrum perturbations, the Kneser-type constants and
the effective masses of zones-type potentials. in: {\it Constructive
theory of functions '84\/}, Sofia 1984, pp. 757--766

\litem{22}
Rofe-Beketov F.S. Kneser constants and effective masses for band potentials.
{\it Sov. Phys. Dokl.\/} {\bf 29} (5) (1984) 391--393

\litem{23}
Schmidt K.M. On the essential spectrum of Dirac operators with spherically
symmetric potentials. {\it Math. Ann.\/} {\bf 297} (1993) 117-131

\litem{24}
Schmidt K.M. Dense point spectrum and absolutely continuous spectrum in
spherically symmetric Dirac operators. {\it Forum Math.\/} {\bf 7} (1995)
459--475

\litem{25}
Schmidt K.M. Critical coupling constants and eigenvalue asymptotics of
perturbed periodic Sturm-Liouville operators. {\it Commun. Math. Phys.\/}
{\bf 211} (2000) 465--485

\litem{26}
Sobolev A.V. Weyl asymptotics for the discrete spectrum of the perturbed
Hill operator. {\it Adv. Sov. Math.\/} {\bf 7} (1991) 159--178

\litem{27}
Weidmann J. {\it Spectral theory of ordinary differential operators.\/}
Lect. Notes in Math. {\bf 1258}, Springer, Berlin 1987

\bye